\documentclass{article}
\usepackage{amsmath, amssymb, graphics, setspace}

\newcommand{\mathsym}[1]{{}}
\newcommand{\unicode}[1]{{}}
\usepackage{graphicx}

\begin{document}

\title{Higher-order root distillers}
\author{M{\' a}rio M. Gra{\c c}a\thanks{Departamento de Matem{\' a}tica ((LAETA/IDMEC),
Instituto Superior T{\' e}cnico, Universidade de Lisboa, Lisboa, Portugal. $mgraca@math.tecnico.ulisboa.pt$.   }
}
 
\maketitle

%\date{}
%M{\' a}rio M. Gra{\c c}a\\
%\\
%Departamento de Matem{\' a}tica\\
%(LAETA/IDMEC)\\
%Instituto Superior T{\' e}cnico\\
%Universidade de Lisboa\\
%Lisboa, Portugal\\
%March 4, 2015.\\
%M{\' a}rio M. Gra{\c c}a
\begin{abstract}
Recursive maps of high order of convergence $m$ (say $m=2^{10}$ or $m=2^{20}$) induce certain monotone step functions from which one can filter relevant information needed to globally separate and compute the real roots of a function on a given interval $[a,b]$. The process is here called a root distiller. A suitable root distiller has a powerful preconditioning effect enabling the computation, on the whole interval, of accurate roots of an high degree polynomial. Taking as model high-degree inexact Chebyshev polynomials and using the {\sl Mathematica} system, worked numerical examples are given detailing our distiller algorithm.
\end{abstract}

\section{Introduction}

By a higher-order root distiller we mean an algorithm to compute, simultaneously, all (or almost all) the real roots of a polynomial $f$
of (high) degree $d$, in a given interval $[a,b]$. The algorithm relies on a single application, on the interval, of a map $g$ of (very) high
order of convergence $m$ (for instance $m= 2^{10}=1024$ or $m= 2^{20}=1\, 048\,576$). The approach may be seen as
a modern computational perspective of the global Lagrange{'}s ideal \cite{La}.\\
\\
Once defined such a higher-order map $g$, the roots of the function $f$ in the interval are obtained from a table $\mathcal{L}= \left\{ (x_1,y_1), (x_2,y_2), \ldots,(x_N,y_N)\right\}$, where $y_i=g (x_i)$, and $x_i$ belong to an uniform grid (of $N+1$ nodes) of width $h$, 
defined on $[a,b]$. For suitable choices of $h$ and $k$  (the parameter $k$ controls the order of $g$), the map $g$ induces an \textit{invariant monotone step function} leading to certain subsets of  $\mathcal{L}$, say $S_1,S_2,\text{...},S_i,\text{...}, S_r$. These subsets will be called \lq{platforms}\rq\ 
(see Section \ref{sec3}) and have the property that the second component of the points $(x,y)_i$ on the platform $S_i$  are equal to (or close to) the  i-th root of the polynomial.
 \\
A graphical inspection of the list $\mathcal{L}$ may be useful not only to observe the distribution of the roots in $[a,b]$ but also to suggest
good choices for the two parameters controlling the algorithm, which are the mesh width $h$ and the index $k$, the latter related to the order of convergence $m$
of the map $g$ ($m=2^{k+1}$ in the case of simple roots). \\
\\
Although in this work we only deal with roots of polynomials, the same procedure can be adapted to non-algebraic equations $f(x)=0$ having
at least one root in a given interval, or with the case of multiple roots or even to functions in $\mathbb{R}^n$ \cite{MG}. 

\medskip
\noindent
In finite arithmetic, one of the main feature of our distiller process is that, by construction, the values $y_i$ are generally quite immune to
rounding error propagation. Therefore a roots{'}s distiller can be seen as a powerful pre-conditioning instrument, in particular for polynomials
whose coefficients are numeric. The referred immunity to rounding error propagation is closely related to the fact that a map of higher-order of
convergence leads necessarily to a stationary {`}monotone machine step function{'} --- if the recursive process which generates the map $g$ is taken
appropriately, that is, for $k$ sufficiently large. Details on the monotone machine step functions are further explained in sections \ref{sec2} and \ref{sec3}.

\medskip
\noindent
Thanks to the super-attracting property of a map of a high order of convergence, each tread (or {`}platform{'}) of the monotone machine step function
-- corresponding to the theoretical subsets \(S_i\) referred above-- contains several machine accurate  values  $y_i$ which are close approximations
of the zeros of the given function $f$. In general, all the necessary information in order to approximate the zeros of a given map with a prescribed
accuracy is contained in these treads.

\medskip
\noindent
Our root distiller is constructed in order to overcome some common numerical issues appearing in the computation of zeros of a given function and
in particular of roots of polynomials of high degree. It is well known the inherent ill conditioning of the computation of polynomial roots. For
instance \textit{Mathematica} commands for approximating roots of polynomials of high degree may produce useless numerical results when low-precision
finite arithmetic is used. On the other hand, dealing with exact polynomials of high degree $d$, say $d\geq 100$, prevents us from using exact arithmetic
due to CPU excessive cost. 
\\
To be more precise, suppose that a numeric expression for the (first kind) Chebyshev polynomial of degree 40 is
defined by the command 
{\tt N[ChebyshevT[40}, {\tt $x],8]$}, where the coefficients are deliberately forced to have 8-digits precision. The commands {\tt Solve},
{\tt Reduce}, and {\tt Roots} produce useless numerical results (cf. paragraph \ref{subsec1}) since the computed roots are heavily contaminated by rounding error (even though
the degree $d=40$ of such polynomial is moderate). Our root distiller deals efficiently not only with this case but it also produces accurate answers,
for instance with a $500$\--degree Chebyshev polynomial.

\medskip
\noindent
In Section \ref{sec2} we detail the construction of a specific map $g$. For that, it is given a positive integer $prec$ and two parameters $h$ and $k$. The
parameter $k$ controls the order of the map $g$ to be constructed, $h$ is the mesh size and $prec$ fixes the precision to be used in the computations of
the images $y_i=g(x_i)$ in the list $\mathcal{L}$.

\medskip
\noindent
In Section \ref{sec3} it is illustrated how a map $g$ of high order of convergence leads to a monotone step function which contains the relevant information
to be distilled. We chose as basic model a 4-degree Chebyshev polynomial of the first kind with $prec=8$ and the parameters $h=0.1$ and $ k=3$. The
respective map $g$ has order of convergence 16 and the absolute error of computed roots is of order $10^{-8}$, meaning that the accuracy used
on $\mathcal{L}$ is preserved. 

\medskip
\noindent
{\sl Mathematica} code is presented in sections \ref{sec2} and \ref{sec4}, including the process used  for filtering the relevant values in the respective list
$\mathcal{L}$ (other filtering possibilities may also be considered). 

\medskip
\noindent
Numerical examples have shown the efficiency of the proposed distillers. In particular, we  construct here a distiller for the computation of the positive
roots of the Chebyshev polynomial of degree 500, defined in $[-1,1]$, with  precision forced to be $prec=5000$ and  parameters  $h=0.00025$ and $k=20$ (and so the respective map $g$ has order of convergence $m=2^{21}= 2 \,097 \,152$). The computed roots have 5000-correct digits.

\medskip
\noindent
An automatic  choice of appropriate parameters $h$ and $k$ in order to achieve a preassigned tolerance error can be done, but this is out of the
scope of the present work.

%==================
\subsection{Motivation: a low precision Chebyshev polynomial}\label{subsec1}

Setting the precision $ prec=8$,  we obtain the following \textit{Mathematica} expression for the Chebyshev{'}s polynomial of degree
40, $N[ChebyshevT[40\,,\,x],\,prec]$:

$$
\begin{array}{l}
 1.0000000-800.00000 x^2+106400.00 x^4-5.6179200\times 10^6 x^6+\\
 +1.5690048\times 10^8 x^8-2.6777682\times10^9 x^{10}+3.0429184\times 10^{10} x^{12}-\\
-2.4343347\times 10^{11} x^{14}+1.4240858\times 10^{12} x^{16}-6.2548083\times 10^{12} x^{18}+\\
+2.1002988\times10^{13} x^{20}-5.4553215\times 10^{13} x^{22}+1.1029237\times 10^{14} x^{24}-\\
-1.7375290\times 10^{14} x^{26}+2.1236466\times 10^{14} x^{28}-1.9918340\times10^{14} x^{30}+\\
+1.4055280\times 10^{14} x^{32}-7.2155451\times 10^{13} x^{34}+2.5426206\times 10^{13} x^{36}-\\
-5.4975581\times 10^{12} x^{38}+5.4975581\times10^{11} x^{40}.
\end{array}
$$

The \textit{ Mathematica} commands {\tt Solve}, {\tt Reduce} and {\tt Roots} produce, respectively, the following useless output:

\medskip
\noindent
 \centerline{\includegraphics[scale=0.95]{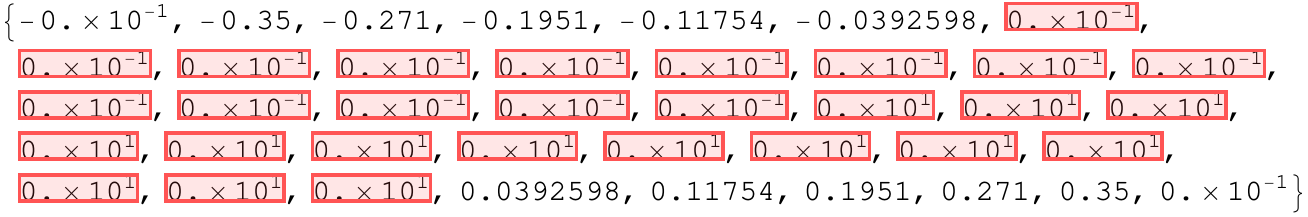}}

\medskip
\noindent
 \centerline{\includegraphics[scale=0.85]{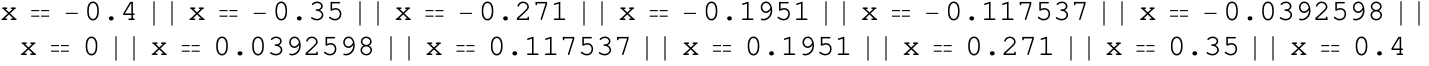} }

\medskip
\noindent \centerline{\includegraphics[scale=0.89]{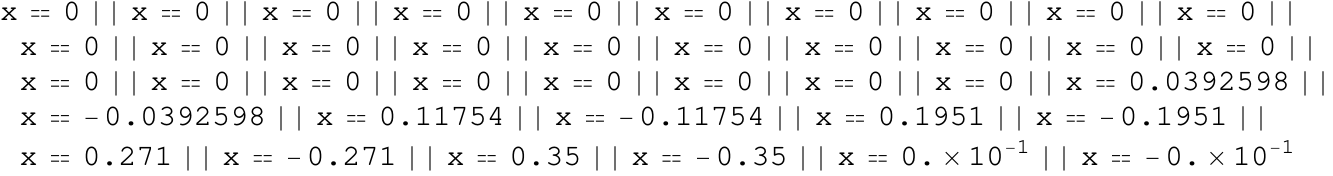}}
 
\noindent Of course the commands {\tt NSolve} and {\tt NRoots} also give useless values.

\medskip
\noindent
We aim to obtain \lq{machine}\rq\  acceptable answers, that is to compute the real roots of the 40-degree Chebyshev polynomial in the interval $[a,b]=[-1,1]$,
which are simple and distinct, with an accuracy close to that of the data (recall that 8-digits precision has been assigned to the polynomial coefficients). 

\medskip
\noindent
The algorithm which
we call \lq{root distiller}\rq\ is described in what follows, and shows to be able to accomplish such a desideratum. The code can easily be included
in a single function in order to produce the referred machine point list $\mathcal{L}$, once predefined the function f, the bounds of the interval, the preassigned
precision $prec$ and the mesh size $h$. After a convenient filtration of the data in $\mathcal{L}$, the respective output should be considered global
in the sense that  it is able to (simultaneously) produce accurate approximations of the roots in the interval, as well as realistic error estimates to
each of them (see Section 4).

%==================================
\section{Higher-order educated maps and monotone step functions}\label{sec2}

In a global approach to roots{'}s computation by means of a smooth  high order
of convergence map $g: [a,b] \subset \mathbb{R}\longrightarrow \mathbb{R}$, many of the domain points are irrelevant, in the sense that their image under $g$ might not be a number or is repealed from a fixed point of $g$. In fact any map of order of convergence greater than one enjoys such a repealing/attracting property -- like
in the well known cases of the Newton's or secant methods for approximation of simple roots. So, once defined a map $g$ of sufficiently high order of
convergence, the points which are not a number, nor in the interval $[a,b]$, neither attracted to a fixed point will be ignored. This is the reason
why we then will call $g$ an \lq{educated}\rq\  higher-order map.

\medskip
\noindent
In general, the recursive maps to be considered have order $m$ of convergence which can go up to $m= 2^{10}=1024$, or greater. The recursive
process used to define the map makes possible to obtain a monotone step function defined in the interval $[a,b]$. From this step function one extracts
the relevant computed $g$-images through a filtering process in order to obtain as output most, or all, the roots of the equation $f(x)=0$. In particular,
our distillers will allow us to compute the roots of a Chebyshev{'}s polynomial of high degree, a task not feasible by the exact methods provided
by the \textit{Mathematica }system (version 10.02.0 running on a Mac OS X personal computer has been used in this work), unless the interval is small and the working precision high.

\subsection{Recursive construction of the higher-order map g}\label{subsec21}

We now explain  the   recursive construction of a map $g$ of high order of convergence by taking as a seed the Newton's map. The map $g$ is the $k$-fold composition
of this seed and has order of convergent $m=2^{k+1}$. The construction of $g$ goes through and \lq{education}\rq \ process aiming  to obtain a map satisfying
a fixed point theorem in the interval [a,b]. More precisely, $g$ is constructed in order to satisfy the following properties: 
\begin{itemize}
\item[(i)] $g([a,b])\subseteq  [a,b]$ .
\item[(ii)] $|g(x)-x|\leqq  (b-a)$, { }for all $x$ for which \ {\tt NumericQ[g(x)]} is \textit{True }. 
\item[(iii)] The points $x\in [a,b]$ not satisfying (i) and (ii) are ignored (a \textit{Null} is assigned to $g(x)$) .
\end{itemize}
\medskip
\noindent
For $k$ sufficiently high, the educated map $g$ will act on $[a,b]$ as a kind of a 'magnet'\  having both good theoretical and computational  properties.
This \lq{magnetic}\rq \ property is better perceived by inspecting a plot of the respective induced monotone step function. 

\medskip
\noindent
It can be proved that for a fixed mesh size $h$, an educated map induces an invariant (or stationary) monotone step function, whenever the
folding parameter $k$ is sufficiently large. This invariant step function will be called the \textit{machine step function} associated to a map
$ g=g[h,k]$.

\medskip
\noindent
Choosing suitable values for $k$, the second component of points on the treads or \lq{platforms}\rq\  of the associated step function contain (by construction)  accurate approximations of the roots of $f(x)=0$. Moreover, the platforms of such step function are automatically sorted in increasing order of their heights, defining so a monotone step function in the interval $[a,b]$. The later filtering process of the data of this step function will hopely solve the referred global Lagrangian root's problem.

\medskip
\noindent
In the following illustrative example an uniform mesh of points, of width $h$, is defined on the domain range $[a,b]$. The {\tt ListPlot} command is used in
order to observe the behaviour of an higher-order educated map $g$ on the referred mesh. 

\medskip
\noindent
Although in this work  the seed used  in  the recursive process is  the Newton's map, any other method of order of convergence greater
than one could be used. For instance, the secant method, of order $(1+\sqrt{5})/2$, and Ostrowsky's methods, of orders $\geq 3$, are other obvious
options.

%================================
\subsection{The map g from the Newton{'}s seed}\label{subsec21}

Fixing a precision $prec$, assume that the numeric expression for a given function $f $ is in memory as well as the bounds $a$ and $b$ of the interval where
the roots of $f(x)=0$ are required. Given the (folding) parameter $k$, the following code defines a general recursive function {\tt g[x,prec]} using Newton{'}s
map as seed (see below functions {\tt newt[0,x,prec]} and its recursive version {\tt newt[k,x,prec]}). The map $g$, of order of convergence $m=2^{k+1}$,
is given below as the function named  {\tt g[x,prec]}.
\\
Note that when $y$={\tt g[x,prec]} is a number, the assigned precision to $y$ is forced to be the same as the precision of $x$. This prevents the \textit{Mathematica} system to correct the output of each calculation of $y$  in the case it occurs of a loss of significant digits. 
The code for the function {\tt g[x,prec]} follows.

\medskip
\noindent

\noindent \centerline{
\includegraphics[scale=0.88]{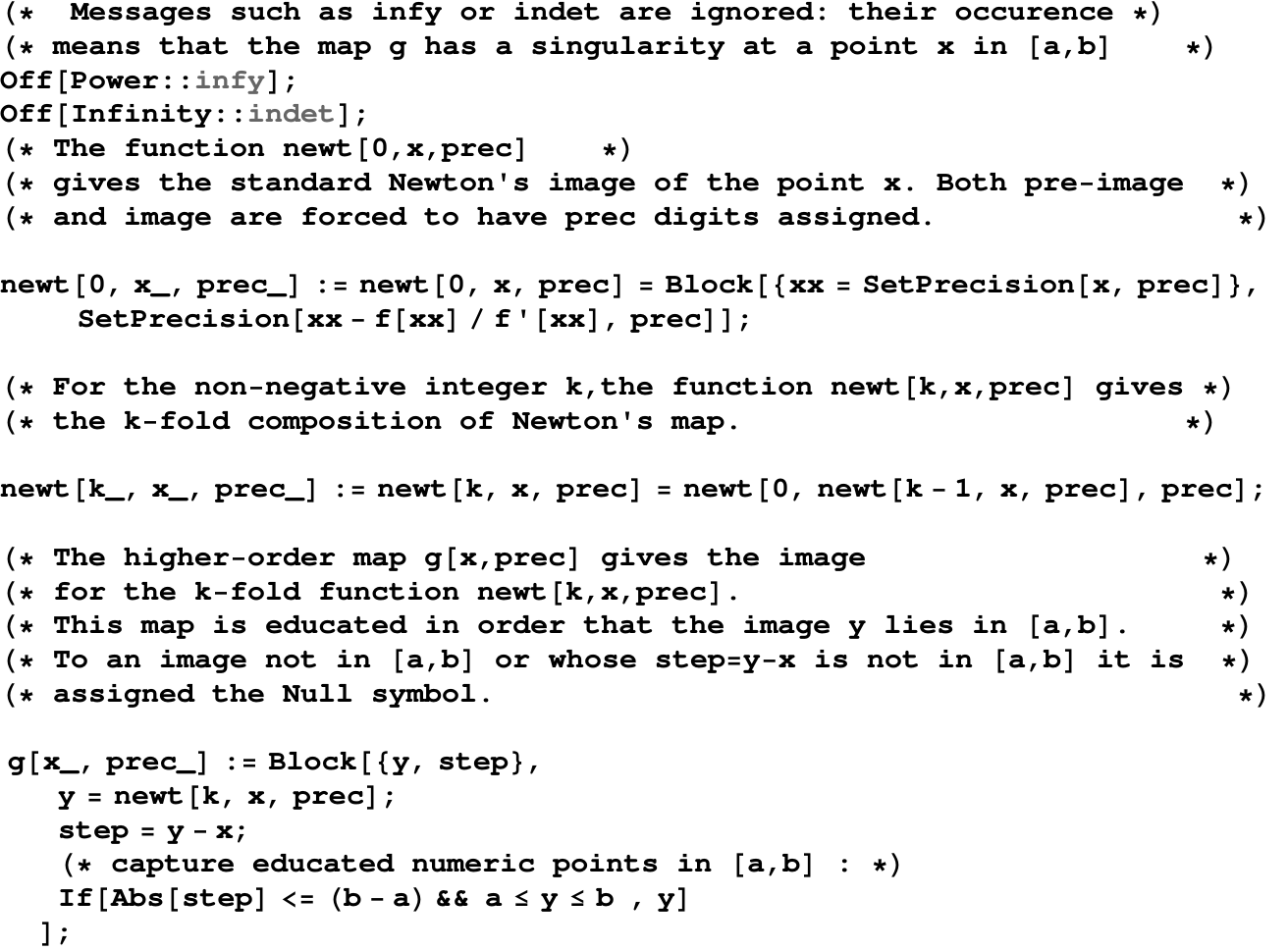}
}
 
 %=============================
\section{An illustration with a low precision 4-degree Chebyshev polynomial}\label{sec3}

As an illustration of the occurrence of a monotone step function induced  by the  map $g$, let us consider a 4-degree Chebyshev polynomial in the interval $[a,b]=[-1,1]$, and  define an uniformly spaced mesh of width $h=0.1$ (that is 21 equally spaced nodes).
We  apply the above {\tt g[x,prec]} code, with assigned parameters $prec=8$ and  $k=3$, that is, in this case the map $g$ has order of convergence $m =2^{k+1}=16$.  

\begin{figure}[hbt] 
\begin{center}
\includegraphics[scale=0.5]{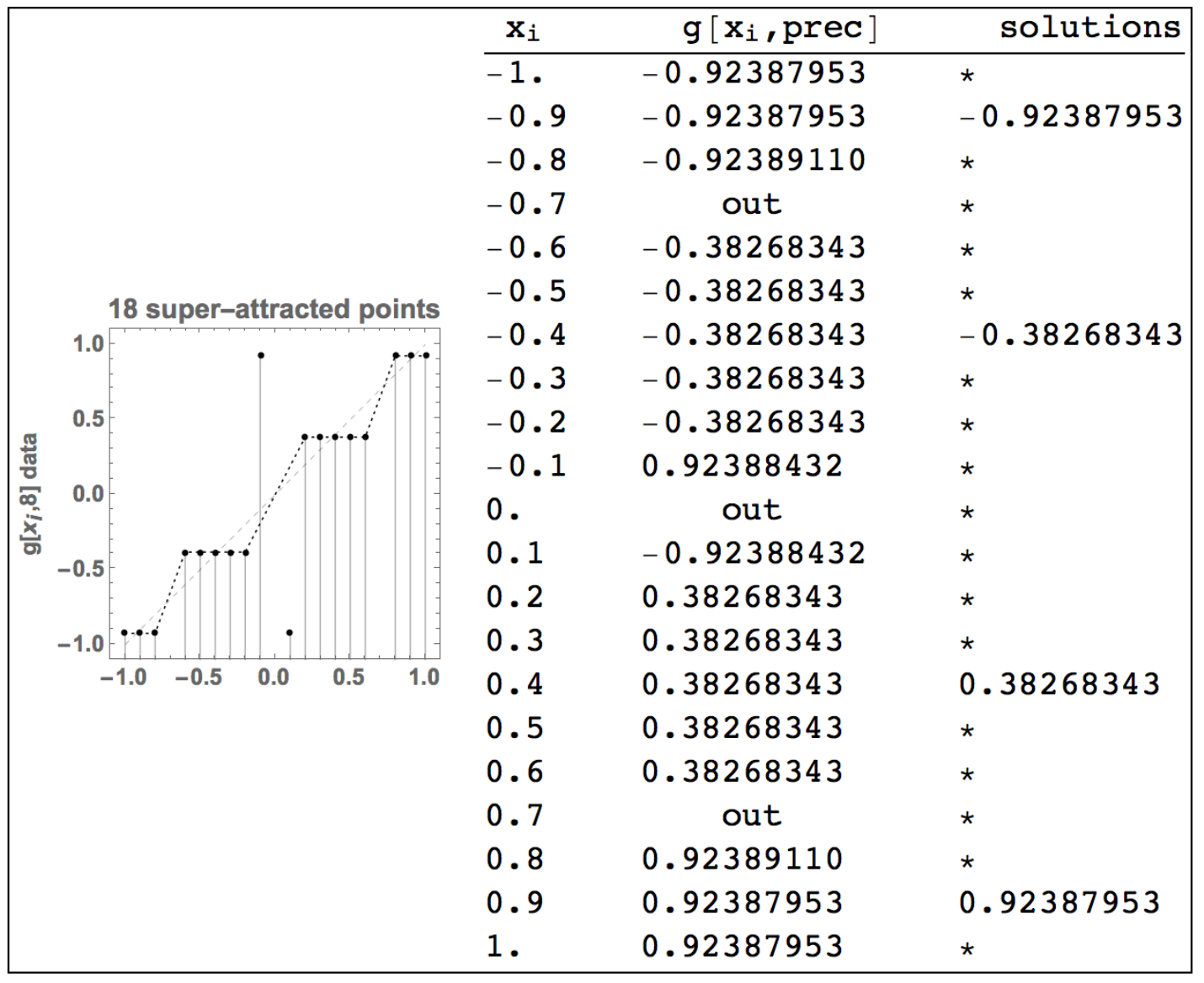}
\caption{A step function associated to a 8\--digits precision, 3\--fold map $g$, applied to a 4\--degreed Chebyshev polynomial. \label{figure1}}
\end{center}
 \end{figure}

\medskip
\noindent
In Figure \ref{figure1} it is displayed the plot of the respective point list $\mathcal{L}$ and  the  computed values $g[x,prec]$. This figure is self-explanatory:  there are 4 roots corresponding to the 4 platforms in the displayed graphic;  an increasing step function is suggested
by the dotted broken line. 
Each of the four observed platforms is formed by a set of points and  from this set of  pairs $(x,y)$ we  filter  the value $y$ of the pair   which has  the  second component  closer to the first. The value of $y$ filtered this way  is an approximation of a polynomial's root. In this example, the 4 points filtered  are displayed in the last column of the table. These 4 values are 8-digit accurate roots of the Chebyshev
polynomial of degree 4.

 \medskip
 \noindent
Looking at the second column of the table in Figure 1,  it is clear  the super-attracting property of the \lq{magnet}\rq\  $g$:  all the numeric points in each
platform have the second component very close to the respective exact  fixed point of the map $g$. 
Repetting the computations for $k=4$, the corresponding table is identical to the one in Figure 1, except  the image of $x= - 0.8$, which is $-0.92387953$.
This means that the respective map $g$, now of order $2^{k+1}=32$, is invariant and so the former computed roots have indeed 8 correct digits.  Saying it  in other words \---  the 4 computed values are \lq{machine\rq\,   fixed points for this map (see the last column in Figure 1).

%================================
\section{Distilling Chebyshev polynomials of high degree}\label{sec4}

For a given precision $prec$ and appropriate parameters $h$ and $k$, we developed a simple \textit{Mathematica} code which will be tested in order to
 approximate the roots of high degree Chebyshev polynomials. A realistic estimative of  the error of each computed root is also easily obtained.
In fact, if $g$ has order of convergence $m>1$ and $y=g(x)$ is a value close to a fixed point $\alpha $, the error of $y$ satisfies $\alpha -y\simeq
 g(y)-y$.

\medskip
\noindent
Concerning our polynomial models, since the roots become closer when we increase the degree of the Chebyshev{'}s polynomial, the value of the
mesh size $h$ and the parameter precision  $prec$ need to be adjusted accordingly. Our aim is to obtain a highly accurate bound for the positive
roots of a 500-degree Chebyshev polynomial of the first kind.

%============================
\subsection{Example: a 500-degree Chebyshev polynomial}\label{subsec41}

A 500-degree Chebyshev{'}s polynomial is highly oscillating and so its roots are very close. Therefore, it is necessary to set a sufficiently
large value of the precision $prec$, and choose a convenient mesh size $h$ in order to obtain the required  numerical results. We consider now the interval to be $[0,1]$.

\medskip
\noindent
 In order to observe the platforms of our distiller, we display the plots of the educated maps $g$ for some values of the parameter $k$ and of $prec$ as a guide for the choice of the right values of these parameters. The next figure compares the graphics of the function
 {\tt g[x, 100]}, respectively for $k=0$ and $k=10$ (the mesh size in the plot is not uniform since it is automatically generated by function {\tt Plot} in the interval $[a,b]$).\\
 \\
It is clear from Figure \ref{figure2} that a $100$\--digits precision is not enough to obtain the roots localised on the right\--half domain, while one can expect a good separation of roots in the interval by setting $prec=5000$, as suggested by Figure \ref{figure3}. 
\begin{figure}[hbt] 
\begin{center}
\includegraphics[scale=0.45]{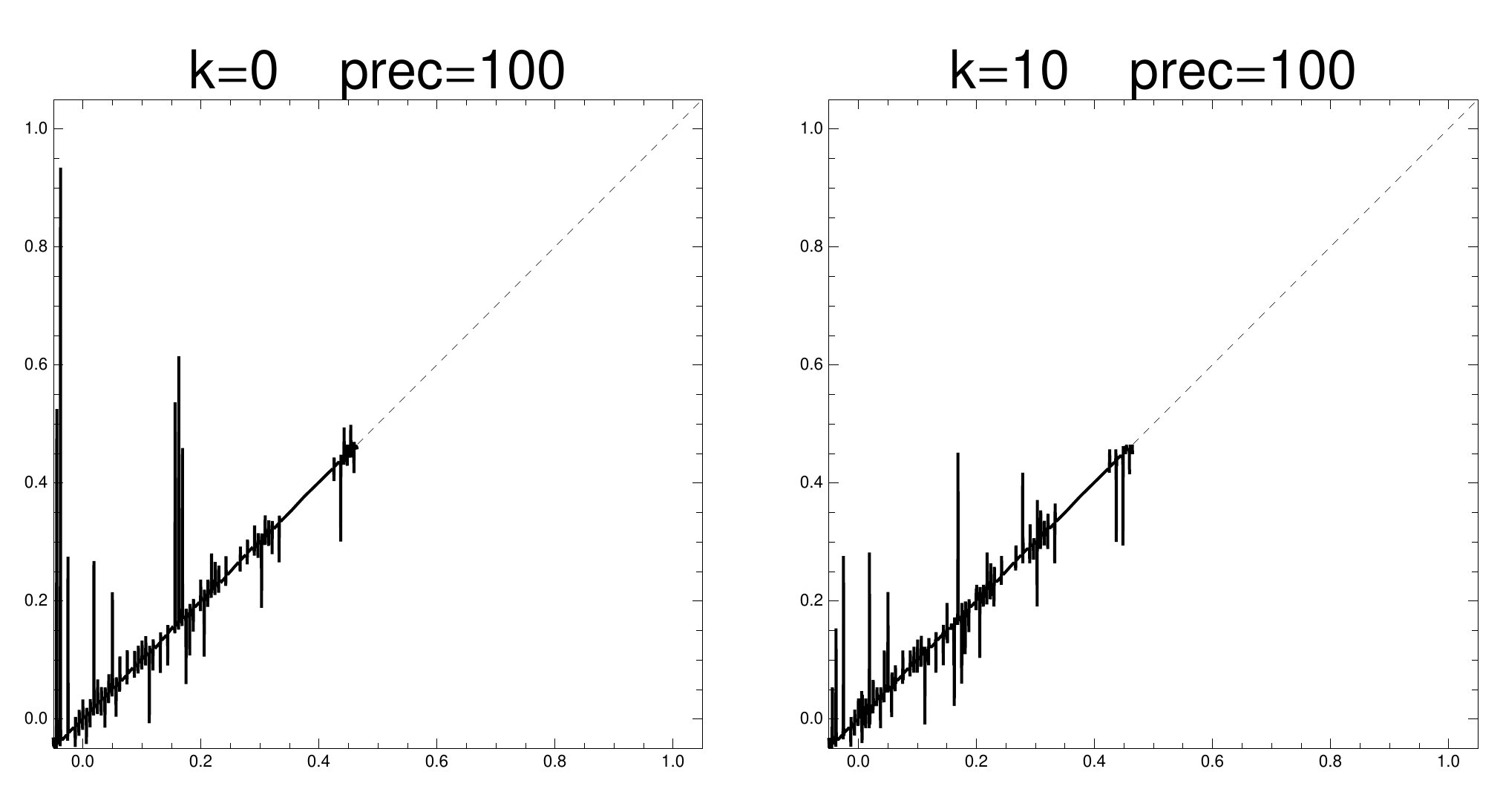} \caption{For $prec=100$ and a mesh automatically generated by $Plot$, only the roots localised on the left half\--domain are detected. \label{figure2}}
\end{center}
 \end{figure}
\\

 \begin{figure}[hbt] 
\begin{center}
\includegraphics[scale=0.45]{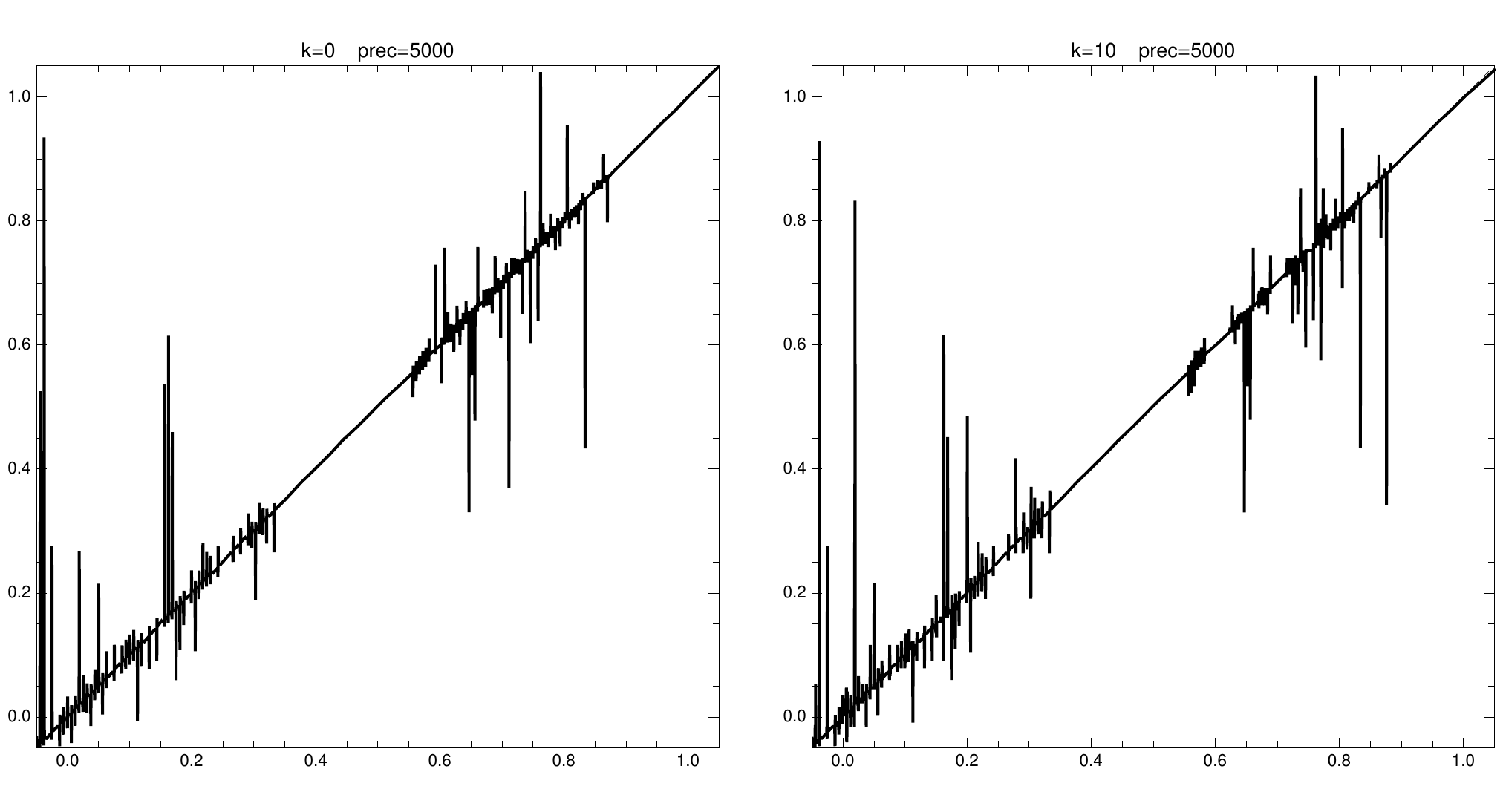}
\caption{Increasing $prec$ to $prec=5000$ and for $k=10$\--folding, the roots on $[0,1]$ of the 500\--degree Chebyshev's polynomial are detected. \label{figure3}}
\end{center}
 \end{figure}

\noindent
An advantage of choosing an higher order of convergence map $g$  will become more apparent if we restrict the interval to $[0.9,1.0]$, which contains
 the desired greatest root of the polynomial. We proceed with the computation of a $5000$ correct digits bound for the positive roots of the polynomial by using a distiller whose folding parameter is $k=20$ (the respective educated map has order $m = 2^{21}=2\,097\,152$). We note that this time the graphics are quickly produced  since we use {\tt ListPlot} instead of {\tt Plot} and therefore  the respective mesh has now much less points than in the usage of  {\tt Plot}.

%====================================================== 
\subsection{ The filtering stage}\label{subsec4}

Efficient \textit{Mathematica} commands for manipulating lists, such as $Cases$, $Parti\-tion$ and $Union$, are particularly useful in the filtering stage
of the distiller algorithm.
\\
Assume that all the numeric points in $\mathcal{L}$ have been assigned to a list named $data$. The first step in the filtration deals with the choice
of points $(x,y)$  sufficiently close to the bisector line, that is  to the line  $y=x$.  We test the condition say, $\left| y-x\right |^2<0.1$,
on the $data$ list, and assign the captured points to a sublist named {\it data1}, as follows

\medskip
{\tt data1} = {\tt Cases[data, x, y; Abs}$[y - x]^2 < 0.1$ {\tt ]; }    (* distille close to bissector *)

\medskip
\noindent
The points in {\it data1} belonging to a certain platform (that is to a \lq{horizontal}\rq\  segment crossing the bisector line) are good candidates for  an
 approximation of the root. The next  sublist, named $data2$, keeps the interesting points. These points are chosen in order to satisfy the condition $(y_1-x_1)(y_2-x_2)<0$, which assures 
 that at least a machine fixed point exists in the interval denoted by $[ x_1, x_2]$:

\medskip
\noindent{\tt data2}= \\
\indent{\small{\tt Cases[Partition[data1,2,1],\{\{x1\_,y1\_\},\{x2\_,y2\_\} \}/;(y1-x1)*(y2-x2)<0 ]; }}

\medskip
\noindent
For a given tolerance, say $tol=10^{-\text{prec}}$, we are interested in filtering the points in the list $data2$  whose second component $y$ differ
from an amount greater than $tol$. Using the command {\tt Union}, we project the platform ignoring machine-duplicate-numbers, obtaining a sublist named $union$,

\medskip
{\tt union = Union[Map[Last, Flatten[data2, 1] ]];}

\medskip
\noindent
The final step in the filtering process checks for the accuracy of the former captured points. First, one filters the values $y$ in the list $union$
for which $f(y)$ is not greater than say $10^{-5000}$. The result is assigned to a sublist called {\it finalA}. Second,  one filters the values $y$ in the list $finalA$, whose estimated absolute error is less than a tolerance, say $tol=10^{-5000}$ . An error bound for each machine root is also computed. The respective code follows.

\medskip
{\tt finalA = Cases[union, y\_ /; f[y] < }$10^{-5000}${\tt];}

{\tt mapf = Map[}$\{${\tt \#, g[\#, prec] - \#$\}$ \&, finalA]; (* roots and error *)}\\

{\tt tol = $10^{-100}$;}

{\tt final = Cases[mapf,$\{$x\_, error\_$\}$/;Abs[error]< tol]; (* error bound *)}

\medskip
\noindent
Assembling the above filtering stage  and the one given at paragraph \ref{subsec21}, a general function is easily obtainable for the whole distiller
algorithm.

%===========================
\subsection{A $5000$\--correct digits approximation for the bound of the positive roots}\label{subsec51}

We now apply our distiller to compute a bound for the positive roots of the 500-degree Chebyshev polynomial  in $[0.99,1.0]$, with an error not exceeding
$10^{-5000}$. Since a bound for the positive roots of the polynomial is required, only the greatest  computed root in the interval will be displayed as well as  its
estimated error.

\medskip
\noindent
In Figure \ref{figure4} the {\tt ListPlot} of the map $g$ is shown, where here  $g$ is the Newton{'}s educated method ($k=0$), the precision is  $prec=5000$, and the mesh size $h=0.00025$.   After filtration an empty list is obtained, meaning that  this map is useless under the previously described filtering criteria.  So a more
powerful \lq{magnet}\rq\  should be used, that is, one needs to  increase the order of convergence of $g$ by taking a greater value of the folding parameter $k$.

 \begin{figure}[hbt] 
\begin{center}
\includegraphics[scale=0.22]{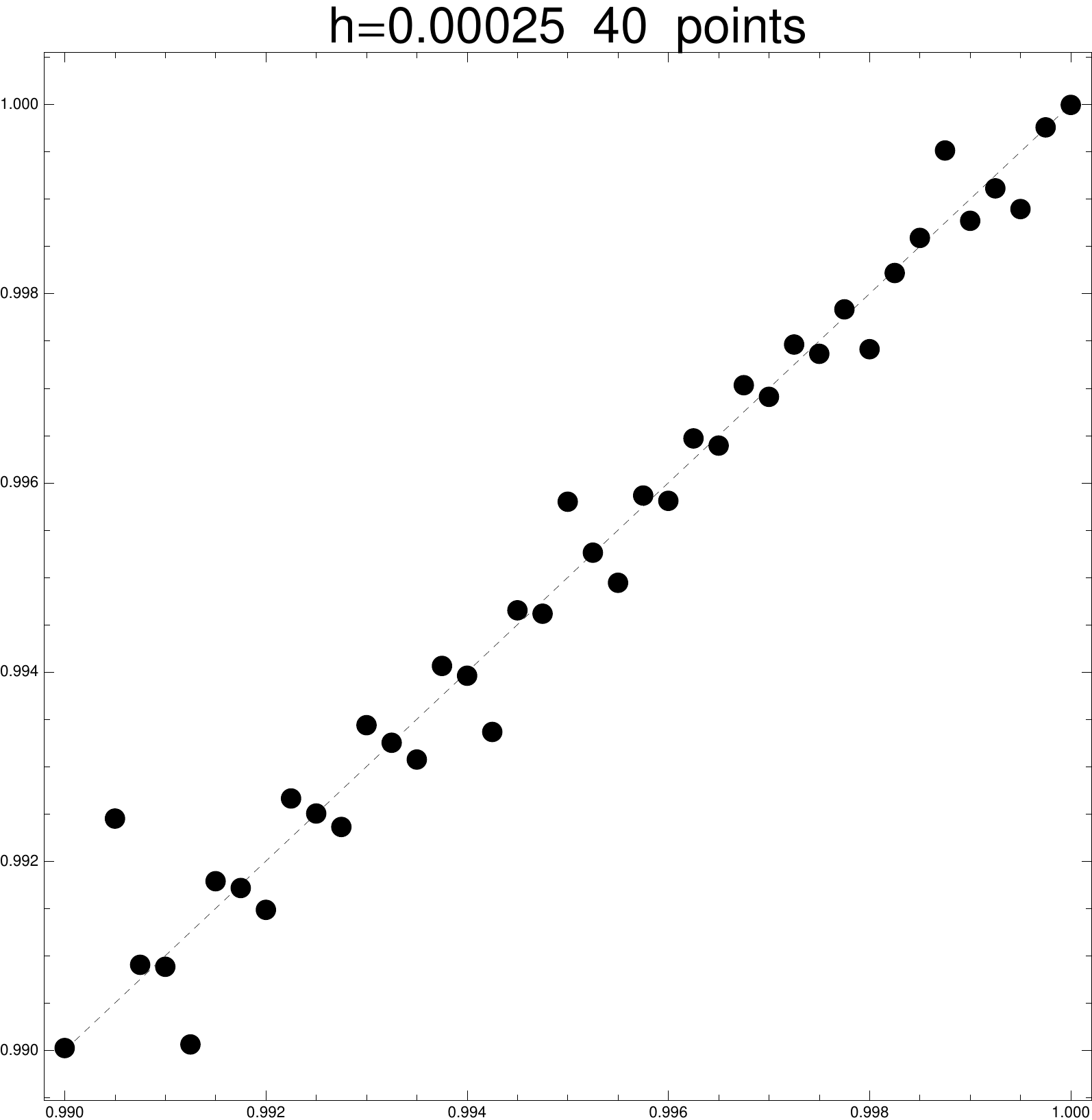} 
\caption{After filtration no points have been captured using the Newton's seed. \label{figure4}}
\end{center}
 \end{figure}

\medskip
\noindent
Increasing to  $k=10$, the respective 2048-order map $g$ (Figure 5 left) enable us to filter relevant points (see Figure \ref{figure5} right) from which high precision
roots can be distilled.  
\begin{figure}[htb] 
\begin{center}
\includegraphics[scale=0.25]{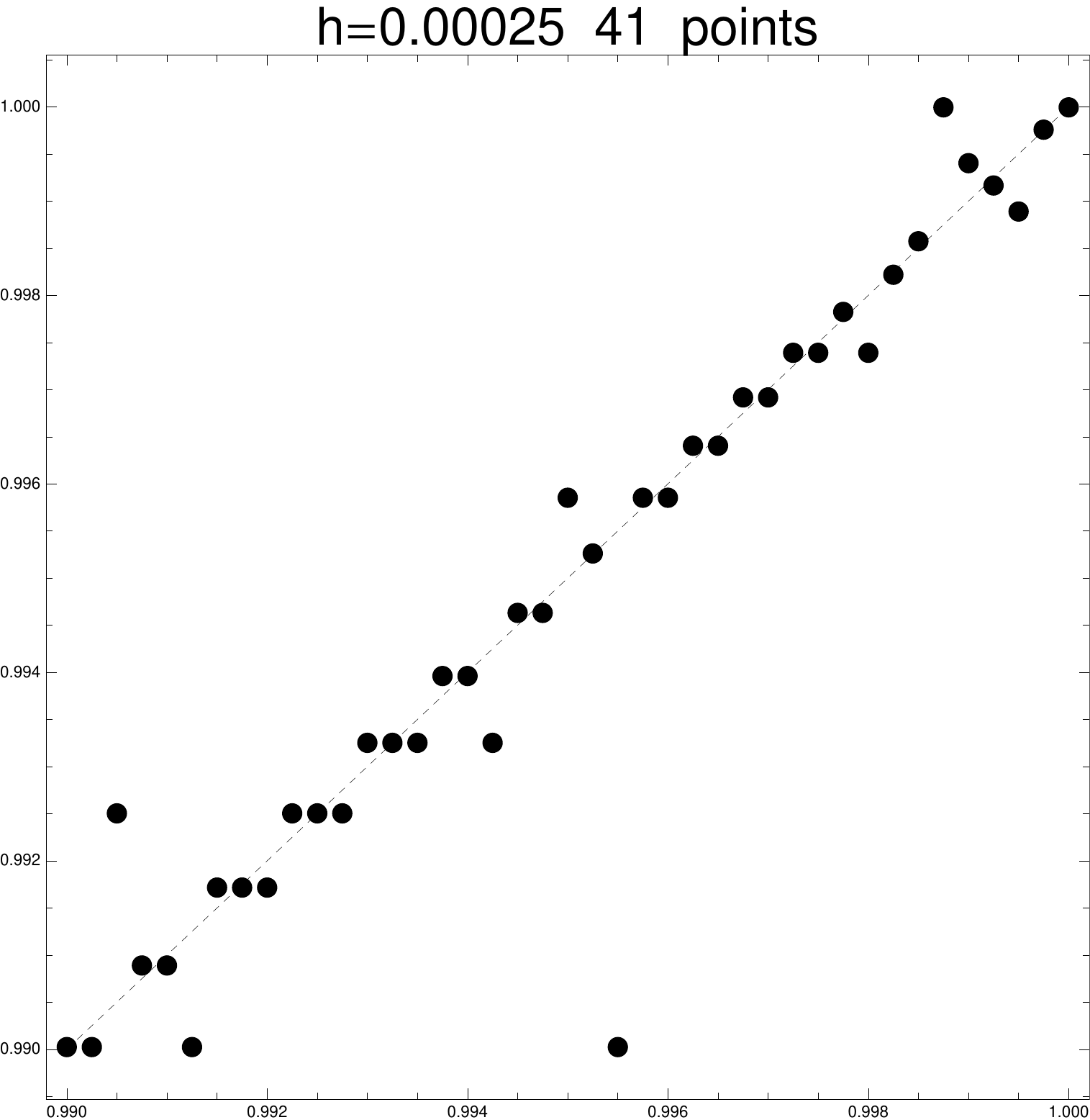} \includegraphics[scale=0.27]{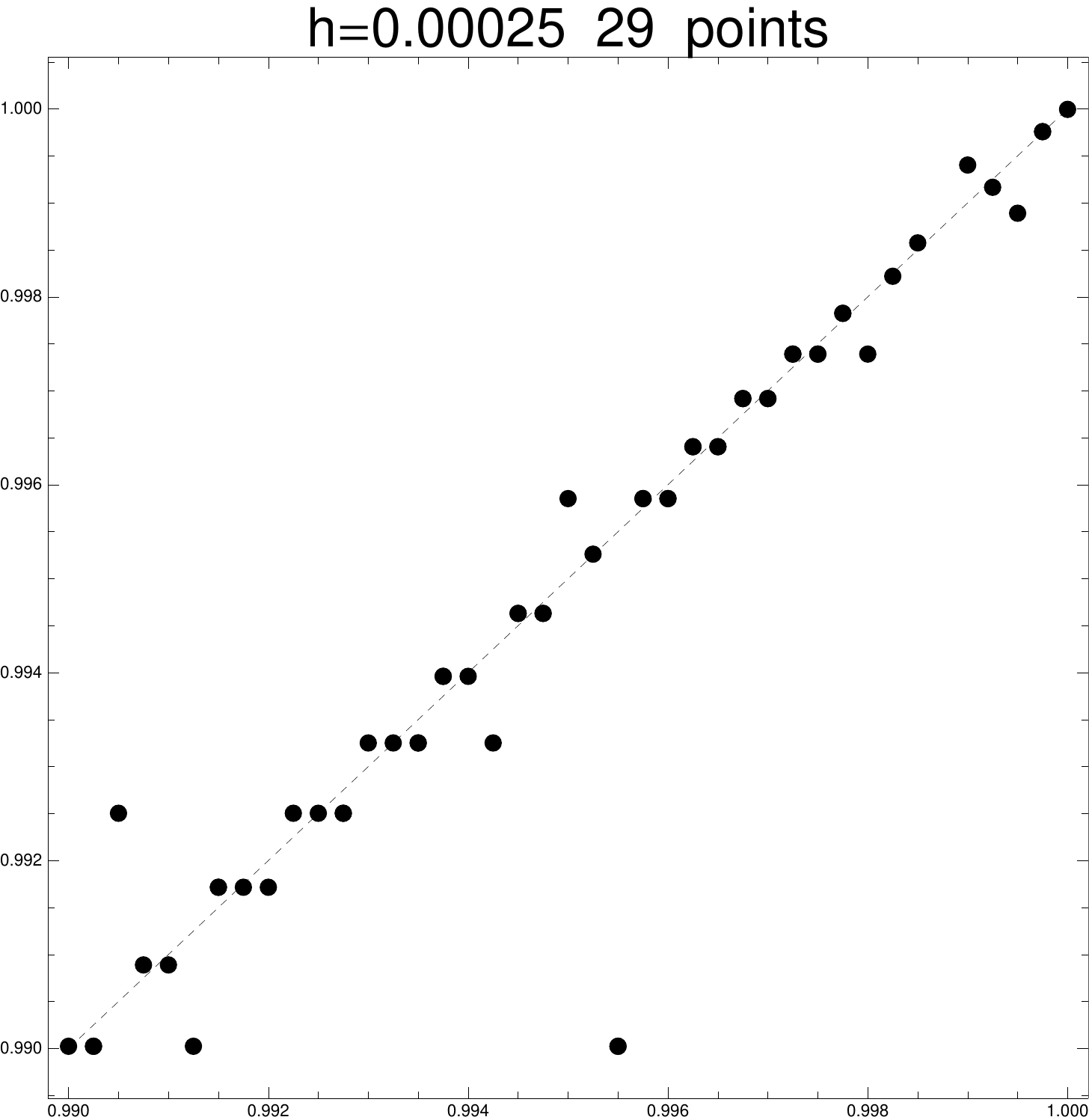}
\caption{For a $10$\--fold map the points in the list  $data2$ of the  respective \lq{machine}\rq\  step function are shown on the right.\label{figure5}}
\end{center}
 \end{figure}
 
\noindent
For $h=0.0025$ and $k=20$ (Figure 6), the complete filtration process leads to 20 fixed points, which are the machine roots of the $500$\--degree Chebyshev polynomial, in the interval $[0.99,1.00]$, for the considered distiller. 

\begin{figure}[htb] 
 \vspace{-0.1cm}
\begin{center}
\includegraphics[scale=0.25]{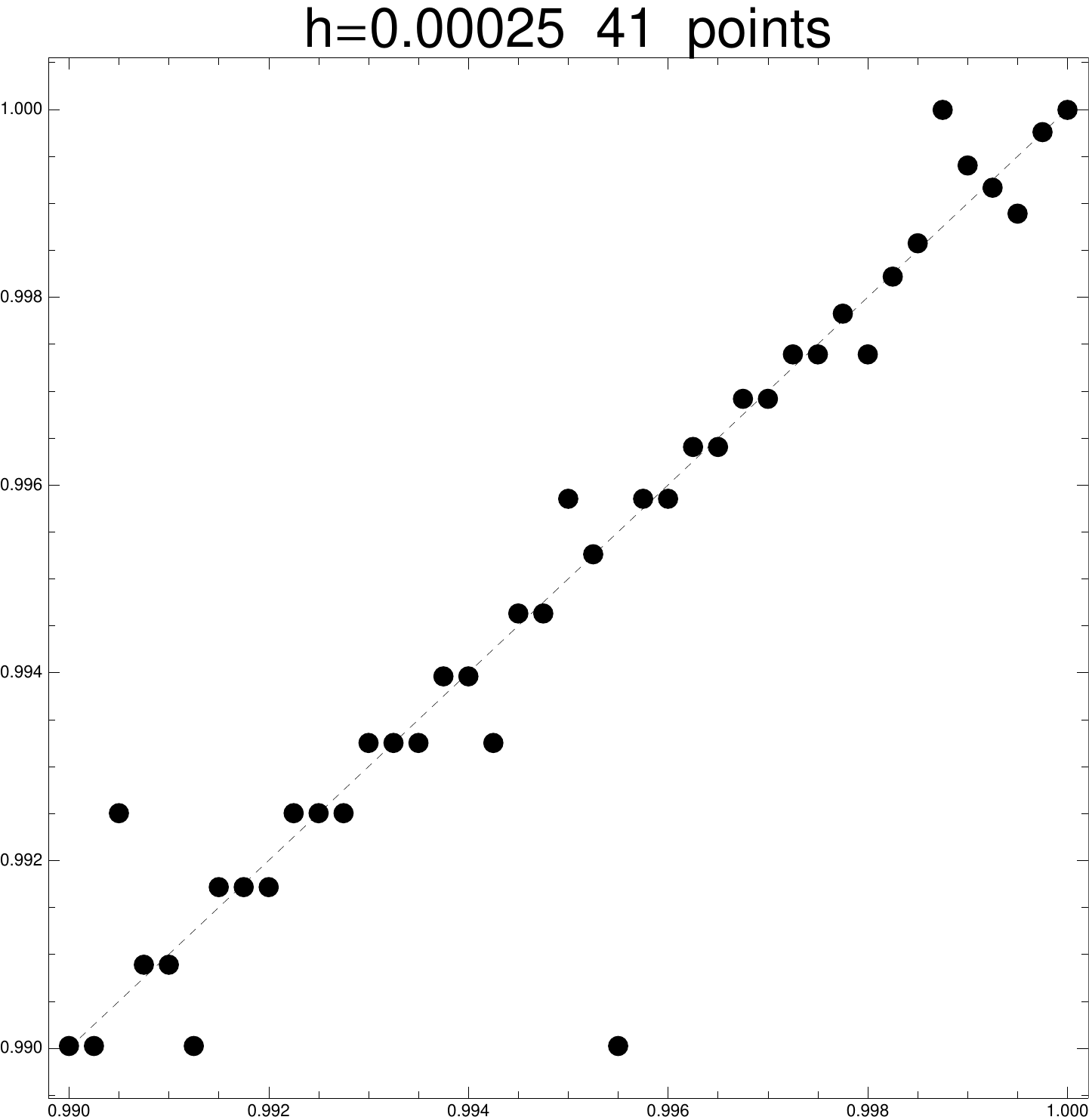}\includegraphics[scale=0.27]{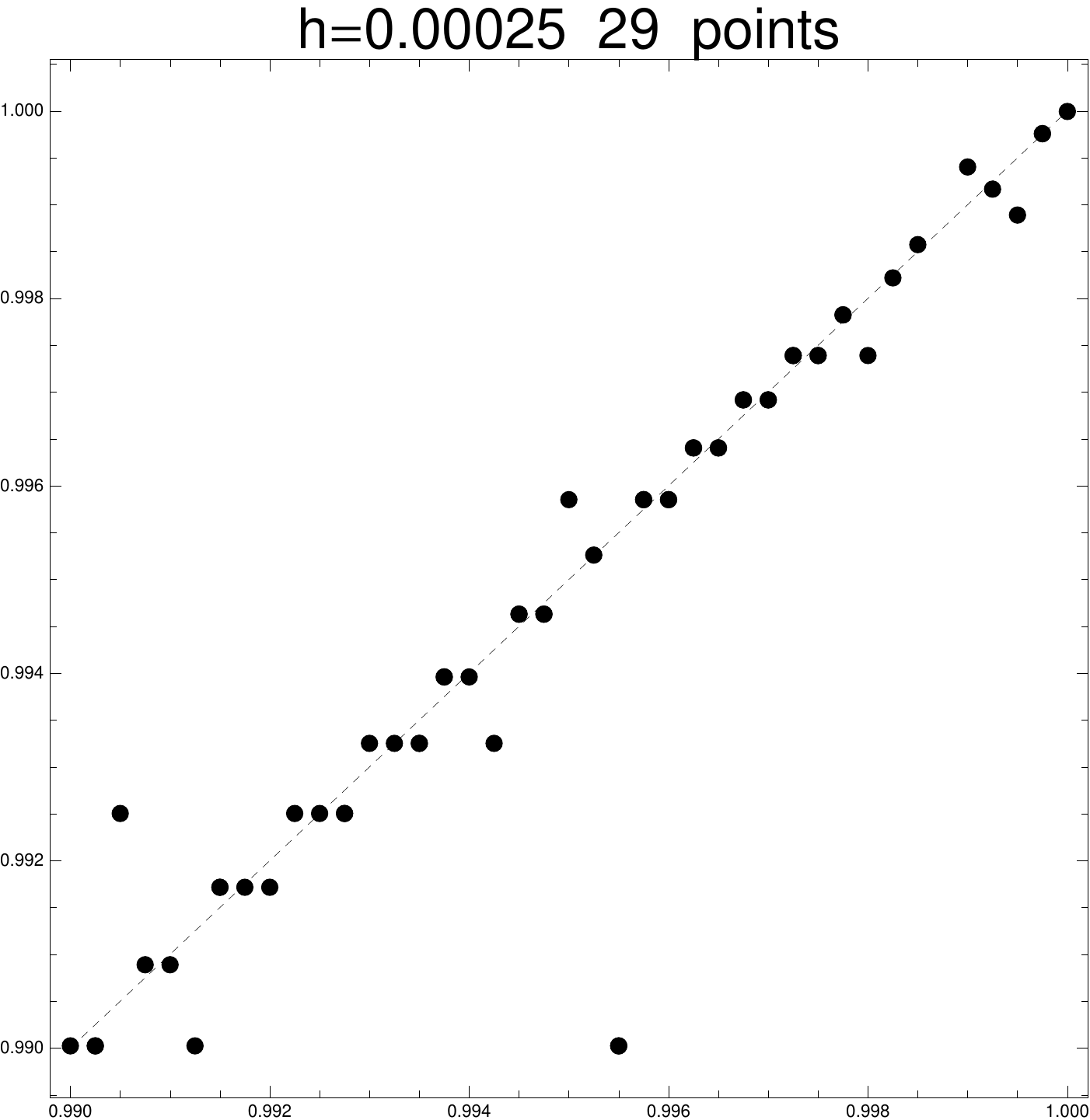}
\caption{Points captured by a 20\--fold $g$ map. \label{figure7}}
\end{center}
 \end{figure}

\medskip
\noindent
Denoting by $\alpha $ the last computed root, a 5000 correct digits bound is obtained. Respectively the first 100 and the last 100 digits of $\alpha
$ are displayed below, as well as its estimated error. 
%\noindent\(\text{}\)
%
%\noindent\(\text{$\alpha $=0.}9999950652018581661118448174487001319149010419592245024005422964693760921636454090733615299842527331\text{...}649619039470741295485455419330860833506167181106321342683094894900371436516812455150525286304558955\)

\begin{center}
\includegraphics[scale=0.55]{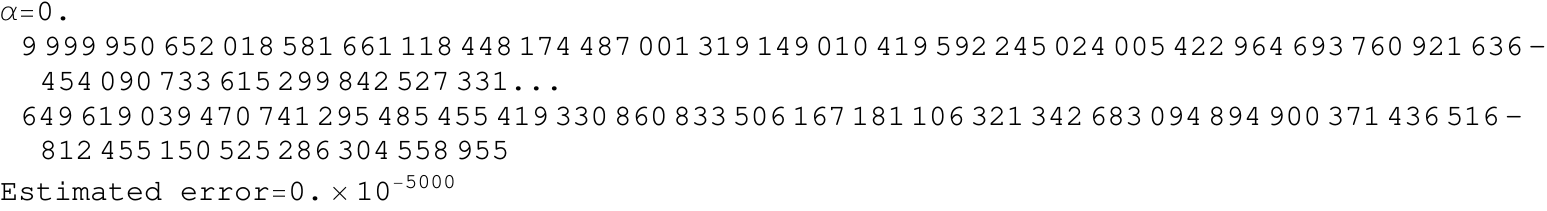}
\end{center}
 Note that a \textit{Mathematica} instruction such as 
$$ x\, /.\, Solve[\{ f[x] == 0, 0.99\leq x\leq 1.0\},x]\}$$
does not produce an answer within an acceptable CPU running time.

 \medskip
 \noindent
Of course, the classical formula giving the zeros of a d-degree Chebyshev polynomial, $\alpha _j=\cos\left( \frac{(2 j+1) \pi }{2 d}\right)$, for $j=0,...(d-1)$,
can be used in order to confirm the above computed value of $\alpha $.

\small{}

\end{document}